\documentclass{article}
\usepackage{amssymb,amsmath,amsthm}

\begin{document}

\large

\theoremstyle{plain}
\newtheorem{theorem}{Theorem}
\newtheorem{lemma}[theorem]{Lemma}
\newtheorem{proposition}[theorem]{Proposition}
\newtheorem{corollary}[theorem]{Corollary}
\newtheorem*{main}{Main Theorem}

\def\mod#1{{\ifmmode\text{\rm\ (mod~$#1$)}
\else\discretionary{}{}{\hbox{ }}\rm(mod~$#1$)\fi}}

\theoremstyle{definition}
\newtheorem*{definition}{Definition}

\theoremstyle{remark}
\newtheorem*{example}{Example}
\newtheorem*{remark}{Remark}
\newtheorem*{remarks}{Remarks}

\newcommand{\sse}{\subseteq}
\newcommand{\nti}{\not \in}
\newcommand{\mn}{\setminus}

\newcommand{\qq}{{\mathbb Q}}
\newcommand{\rr}{{\mathbb R}}
\newcommand{\nn}{{\mathbb N}}
\newcommand{\zz}{{\mathbb Z}}

\newcommand{\al}{\alpha}
\newcommand{\be}{\beta}
\newcommand{\ri}{r_i}
\newcommand{\si}{s_i}
\newcommand{\rp}{\rr^+}
\newcommand{\ep}{\epsilon}
\newcommand{\lam}{\lambda}
\newcommand{\de}{\delta}

\newcommand{\ra}{\rangle}
\newcommand{\la}{\langle}
\newcommand{\sem}{S(}

\title{Positive rational solutions to $x^y = y^{mx}$: \\ a number-theoretic excursion}

\author{Michael A. Bennett \\ University of Illinois at
Urbana-Champaign\\ University of British Columbia\\ \\Bruce Reznick
\\ University of 
Illinois at Urbana-Champaign}

\maketitle

\section{Prologue}

Late in the last millennium, the
second author ran a 
seminar course 
for undergraduates which was intended to introduce them to
problem-solving and question-asking in the context of mathematical
research. He led them through  the classic ``difficult" equation
\begin{equation} \label{eq1}
x^y = y^x,\qquad x, y > 0, 
\end{equation}
whose solution is much easier than one would think at first glance.
Solutions were sought, first over $\rr$, then over $\zz$ and finally
over $\qq$. In the context of this seminar, it was natural to consider
the variant equation
\begin{equation} \label{eq2}
x^y = y^{2x},\qquad x, y > 0, 
\end{equation}
which does not seem to have appeared in the literature. It turns out
that there are solutions to (\ref{eq2}) which do not fit the well-known
parametric pattern of (\ref{eq1}); c.f. (\ref{eq5}) below.  For example, 
\begin{equation} \label{eq3}
x = \left(\frac 45\right)^{128}, \quad y =  \left(\frac 45\right)^{125}
\end{equation}
is a solution to  (\ref{eq2}). This preposterous fact is trivial to
verify: simply substitute (\ref{eq3}) into  (\ref{eq2}), take logs and
transpose:
$$
\frac {2x}y = 2\left(\frac 45\right)^3 = \frac{128}{125} = \frac{\log
x}{\log y}.
$$
As we shall see, the ultimate explanation for this identity is that
$2\cdot 4^3 = 5^3 + 3$. 
Upon discovering (\ref{eq3}), the second author realized he was  in
over his head and 
contacted the first author. This paper is the result.

\section{$x^y = y^x$}

First, we review the familiar, but beautiful solution to
(\ref{eq1}), reserving historical references to the last paragraph
of the section.  We acknowledge the solutions $x=y$ and now let $y = tx$,
$t \neq 1$, so that  
\begin{equation} \label{eq4}
x^{tx} = (tx)^x \implies x^t = tx \implies x =
t^{\frac 1{t-1}}.
\end{equation}
The positive real solutions to (\ref{eq1}) are thus
\begin{equation} \label{eq5}
(x,y) = (u,u) \qquad \text{and} \qquad (x,y) = ( t^{\frac 1{t-1}},
t^{\frac t{t-1}}). 
\end{equation}
(We might equally well have set $y = x^r$  in (\ref{eq1}), and drawn
essentially the same conclusion.)
Since (\ref{eq1}) 
implies that $x^{1/x} = y^{1/y}$ and since $f(u) = 
u^{1/u}$ increases on $(1,e)$ and decreases on $(e,\infty)$, for each
$x \in (1,e)$, there is 
exactly one $y \in (e,\infty)$ so that (\ref{eq1}) holds. In
particular, the only integer solution to  (\ref{eq1}) is $2^4 =
4^2$.

Euler already noted that if $t = 1 + \frac 1n$ for integral $n$, then 
(\ref{eq4}) gives a rational solution to (\ref{eq1}); namely,
\begin{equation} \label{eq6}
x_n = \biggl(1 + \frac 1n\biggr)^n,\quad 
y_n =  \biggl(1 + \frac 1n\biggr)^{n+1}.
\end{equation}
(Observe that, as $n \to \infty$, we have $t \to 1$ and $(x_n)$ and $(y_n)$
increase and decrease monotonically to $e$, as is familiar from
calculus.)
 
To show that these are the only rational solutions, we need an
elementary lemma, whose proof, relying upon the Fundamental Theorem
of Arithmetic, we omit:
\begin{lemma}
Suppose $a, b, m, n$ are integers, with $\gcd(a,b) = \gcd(m,n) =
1$ and $b, n \neq 0$. Then $\left( \frac mn 
\right)^{a/b}$ is rational if and only if $m$ and $n$ are $|b|$-th
powers of integers.
\end{lemma}

Let us now proceed to find all rational solutions to (\ref{eq1}).
By symmetry, we may assume that $y > x$, so  $t >1$. If $x$
and $y$ are rational, then so is $t = y/x$. Write $t$ in lowest terms as
$$
t = \frac pq := \frac {q+d}q,
$$
($d, q > 0$), so that $\frac 1{t-1} = \frac qd$ and  $\frac t{t-1}
= \frac {q+d}d$. With this substitution, we have
$$
x = \left(\frac {q+d}q\right)^{q/d}, \qquad
y =\left(\frac {q+d}q\right)^{q/d + 1}. 
$$
Since $\gcd(d,q) = \gcd(q,q+d) = \gcd(d,q+d) =1$, 
the integers $q$ and $q+d$ must both be $|d| = d$-th powers by  Lemma 1.
 This causes no
problem when $d=1$, of course, and, setting $q = n$, we recover
(\ref{eq6}). Suppose $d > 1$, and write $q = a^d$, $q+d = b^d$ for
positive integers $a < b$, so that $b^d - a^d = d$. Observe that
$$
b^d - a^d \ge (a+1)^d - a^d \ge 1 + da \ge 1 + d > d
$$
for positive integers $a$ and $b$, and so there are no solutions with
$d > 1$.
Finally, we remark that if $n = -r$ is allowed to be negative in
(\ref{eq6}), then $(x_{-r},y_{-r}) = (y_{r-1},x_{r-1})$:
$$
 \biggl(1 - \frac 1r\biggr)^{-r} = \biggl(1 + \frac
 1{r-1}\biggr)^{r},\quad  
 \biggl(1 - \frac 1r\biggr)^{-r+1} = \biggl(1 + \frac 1{r-1}\biggr)^{r-1}.
$$
Thus, by removing the restriction that $n$ be positive, we can
eliminate the constraint $y > x$ in (\ref{eq6}).

According to Dickson \cite[p.687]{D}, the first reference to
(\ref{eq1}) was in a letter from D. Bernoulli to C. Goldbach, dated
June 29, 1728. Bernoulli asserts, without proof, (see
\cite[p.263]{B}) that this 
equation has only one solution in positive integers, and infinitely
many rational solutions. The first person to write about
(\ref{eq1}) in detail was Euler (see \cite[pp.340--341]{E}. Euler made the
substitution $y=tx$ and solved the equation over $\rr_+$ and $\zz_+$,
and presented the rational solutions (\ref{eq6}), without claiming
that they were the only ones. Dickson mentions other writers who went
over the same ground, and then, mysteriously, writes ``*A. Flechsenhaar
and R. Schimmack discussed the rational solutions". The asterisk means
that the paper was ``not available for report"  (\cite[p.xxii]{D}).
These papers appeared in 1911 and 1912 in the journal
\emph{Unterrichtsbl\"atter f\"ur Mathematik und Naturwissenschaften},
and the authors found them without difficulty on the shelves of the
magnificent 
Mathematics Library of the University of Illinois at Urbana-Champaign,
then and now only 200 kilometers south of Dickson's
office. (For graduates of the US school system, a kilometer is a
Canadian mile, attractively priced at 
$0.621$ American miles). In any event, Flechsenhaar \cite{F} appears to
earn credit as the first author to solve (\ref{eq1}) over positive
rationals. In 1967, Hurwitz \cite{Hu} gave the first readily accessible
analysis of the rational solutions; subsequent work on this equation, including
generalizations to algebraic solutions, can be found in \cite{Ha}, \cite{Mit},
\cite{Sa}, \cite{Sa2}, \cite{St} and \cite{Sv}.

\section{$x^y = y^{mx}$}

Let us now consider the generalization of Euler's equation to
\begin{equation} \label{eq7}
x^y = y^{mx}
\end{equation}
where $m>1$ is a fixed positive integer. We will again restrict our
attention to positive solutions $(x,y)$. 
If $x=1$ or $y=1$, then necessarily $(x,y)=(1,1)$.
Supposing that $x, y \neq 1$ and taking logarithms in (\ref{eq7}), we
find that
$$
\frac{m \log y}{y} = \frac{\log x}{x},
$$
so $x \neq y$. Write $y = x^r$, $x\neq 1,\ r \neq 1$. Then $x^{r-1} = m
r$. Therefore, the positive real solutions to  (\ref{eq7}) are given
by $x = y = 1$ and  
$$
x = \left( mr \right)^{\frac{1}{r-1}} \; \mbox{ and } \; y = \left( mr \right)^{\frac{r}{r-1}}.
$$

We now restrict our attention to positive rational solutions.
Since $y/x = m r$, it follows that $r > 0$, and, further, since $x$
and $y$ are rational, we have $r 
\in \mathbb{Q}$.
Let us write $r = \frac ab$ where $a, b \in \mathbb{N}$ with $\gcd (a,
b) = 1$ and set $k = |a-b| \in \mathbb{Z}_+$.  To have $x \in
\mathbb{Q}$, we require  
$$
x = \left(\frac{ma}b \right)^{\frac{b}{a-b}} \in \mathbb{Q}.
$$
Suppose that $\gcd(b,m) = d$, and write $b = db'$ and $m = dm'$.
Then $\gcd(am', b') = 1$ and so we need $am'$ and $b'$ to be $k$-th powers
of integers. If
$$
am' = u^k \; \mbox{ and } \; b' = v^k,
$$
it follows that
\begin{equation} \label{key}
\left| u^k - m v^k \right| = |am' - mb'| = \left| \frac {am}d - \frac
{mb}d \right| = \frac{m}{d}\cdot k = m'k.
\end{equation}

For a given positive integer $m$, we will classify the set $S(m)$ of positive rational solutions to
(\ref{eq7}) as follows:
write
$$
S(m) = \bigcup_{k=0}^{\infty} S_k (m),
$$
where $S_k(m)$ represents the set of solutions $(x,y)$ corresponding
to equation (\ref{key}), 
for $m'$ ranging over all positive integral divisors of $m$. Here,
$S_0(m)$ denotes the solutions 
with $x=y$ (i.e. just the set $(x,y)=(1,1)$ for $m > 1$).

The remainder of this paper will be devoted to analyzing $S(m)$. As a
consequence, 
we will show how this set may be completely characterized for any
given $m$. In case $m=2$ or $3$ 
(where all features of interest for general $m$ may in fact be
observed), we have the following results: 

\begin{theorem} \label{ruski}
If $x$ and $y$ are positive rational numbers for which
$x^y = y^{2x}$ then either 

(a) $x = \left( 2 + \frac{2}{n} \right)^{n}$ and $y = \left( 2 +
\frac{2}{n} \right)^{1 + n}$, and $n \in \zz$, $n \neq 0,-1$;

or

(b) $(x,y) = (1,1), (2, 16)$ or $ \bigl( \bigl( \frac{4}{5}
\bigr)^{128}, \bigl( \frac{4}{5} \bigr)^{125} \bigr)$.
\end{theorem}
and
\begin{theorem} \label{ruski2}
If $x$ and $y$ are positive rational numbers for which
$x^y = y^{3x}$ then either 

(a) $x = \left(3 + \frac{3}{n} \right)^{n}$ and $y = \left( 3 +
\frac{3}{n} \right)^{1 + n}$, and $n \in \zz$, $n \neq 0,-1$;

or

(b) $x = \bigl(\frac {3w_n}{v_n}\bigr)^{v_n^2}$ and 
$y = \bigl(\frac {3w_n}{v_n}\bigr)^{3w_n^2}$, $0 \le n \in \zz$;

or 

(c) $x = \bigl(\frac{w_n}{v_n}\bigr)^{3w_n^2}$ and 
 $y = \bigl(\frac{w_n}{v_n}\bigr)^{v_n^2}$, $0 \le n \in \zz$.

In (ii), $v_n$ and $w_n$ are the integers  defined by 
$
v_n + w_n \sqrt{3} = ( 1 + \sqrt{3} ) ( 2 + \sqrt{3} )^n.
$
\end{theorem}

Our proofs of Theorems \ref{ruski} and \ref{ruski2} (perhaps rather surprisingly)
involve techniques from  
transcendental number theory and Diophantine approximation. 
Moreover, they provide us with an opportunity to illustrate
a fairly diverse grab-bag of methods for solving Diophantine problems.

\section{The cases $k=1$ and $k=2$}

The set $S_1(m)$ is easily computed. Indeed, we immediately find that
$k=1$ implies that either 
$$
x= \left( m + \frac{m}{n} \right)^{n}, \hskip2ex y = \left( m + \frac{m}{n} \right)^{1 + n}
$$
or
$$
x= \left( m - \frac{m}{n} \right)^{- n}, \hskip2ex y = \left( m - \frac{m}{n} \right)^{1 - n}
$$
for $n$ a positive integer (with $n \geq 2$ in the latter case). We note that these solutions
correspond to the  parametrized family (\ref{eq6}) of solutions to
Euler's original equation and provide us with part (a) of Theorems
\ref{ruski} and \ref{ruski2}.

If $k \geq 2$, the situation becomes rather more interesting, though
the set $S_2(m)$ is also not too difficult to describe: it is either empty
or infinite. The following lemmata provide sufficient conditions for
the former to occur. As usual, for $x \in \nn$, let $\nu_2(x)$ be the
largest integer such that 
$2^{\nu_2(x)}$ divides $x$.
\begin{lemma} \label{lem1}
If $\nu_2(m)$ is odd and $k$ is even, then $S_k(m)$ is empty.
\end{lemma}
\begin{proof}
We have $\gcd (a,b) = 1$, and since $k = |a-b|$ is even, $a$ and $b$
must both be odd. It follows that $d = \gcd(b,m)$ is odd as
well. Write, as before, $b = db'$, $m = dm'$, $am' = u^k$, $b' =
v^k$. Then $\nu_2(m) = \nu_2(d) + \nu_2(m')$ is odd, hence so is 
$\nu_2(am') = \nu_2(u^k) = k\nu_2(u)$, which contradicts $k$ being even.
\end{proof}

\begin{lemma} \label{lem2}
If $m = 2^{\nu_2 (m)} \cdot m_1$ is a positive integer for which $m_1 \equiv 1 \mod{4}$, then
$S_2(m)$ is empty. 
\end{lemma}

\begin{proof}
If $\nu_2(m)$ is odd, this is a special case of Lemma \ref{lem1}.  
Since $k=2$ is even, we may conclude as before that $a$ and $b$ are
both odd. Suppose
$\nu_2(m)$ is even, say $\nu_2(m) = 2 t$.
It follows that
$$
u^2 - 2^{2t} m_1 v^2 = \pm 2^{2t+1}\cdot \tfrac{m_1}d
$$
where $u$ and $v$ are coprime and $d$ divides $m$. Since $d\ | \
b$, it follows that $d$ is odd, and so  $d\ |\  m_1$. We thus have
that $2^{t}\ |\ u$, say 
$u = 2^{t} u_1$, whereby
\begin{equation} \label{rufus}
u_1^2 - m_1 v^2 = \pm 2\frac{ m_1}d.
\end{equation}
Since the right hand side of this equation is even, $u_1$ and $v$ have the same parity, and, since $u$
and $v$ are coprime, are both necessarily odd. Since $m_1 \equiv 1 \mod{4}$, this implies that the left
hand side of (\ref{rufus}) is divisible by $4$. Since the right hand side of this equation is congruent to $2$
modulo $4$, this yields the desired contradiction.
\end{proof}

It is a well-known fact 
that if $m$ is a positive nonsquare integer and $c$
is a nonzero integer, then a single solution in positive integers $x$
and $y$ to the equation $x^2- m y^2 = c$ implies the existence of
infinitely many
such solutions. In fact, one can find a finite collection of
pairs of positive integers, say
$$
(x_1,y_1), \; (x_2,y_2), \ldots, (x_r,y_r),
$$
with $x_i^2-my_i^2=c$ for $i = 1, 2, \ldots, r$, such that {\it{every}}
solution in positive integers $(x,y)$ to the equation $x^2-my^2=c$
satisfies
\begin{equation} \label{goofy}
x+ y \sqrt{m} = \left( x_i + y_i \sqrt{m} \right) \cdot \left( u_1 + v_1 \sqrt{m} \right)^k
\end{equation}
where $k$ is a nonnegative integer, $i \in \{1, 2, \ldots, r \}$ and $(u_1,v_1)$ is the smallest
positive integer solution to the equation $u^2-mv^2=1$. The integer $r$ here depends upon $c$ and, potentially,
upon $m$.

It follows that if $S_2(m)$ is not empty, then it is
infinite.  From the theory of Pell equations (or Fermat-Pell equations if
one likes) these correspond to elements of finitely many recurrence
sequences (see Nagell \cite{Na} for a nice, affordable exposition 
of such matters). This fact is also a pretty easy consequence of (\ref{goofy}).
We note that consideration of the case $m=39$ 
demonstrates that the above lemmata do not in fact provide necessary conditions for $S_2(m)$ to be empty.
A routine computation shows that for
$2 \leq m \leq 50$, $S_2(m)$ is  infinite precisely for
$$
m \in \{ 3, 7, 11, 12, 15, 19, 23, 27, 28, 31, 35, 43, 44, 47, 48 \}.
$$

It is also not hard to provide sufficient conditions for
$S_2(m)$ to be nonempty. An  
example of such a result is as follows:
\begin{lemma} \label{lem3} 
If $p \equiv 3 \mod{4}$ is prime then $S_2(p)$ is infinite.
\end{lemma}
\begin{proof}
From a venerable (after citing Euler, we can hardly call this old!)
result of Petr \cite{Pe}, precisely one of
$$
x^2-py^2=-2 \; \; \mbox{ or } \; \; x^2-py^2=2
$$
is solvable in integers $x$ and $y$. This implies, in either case,
that $S_2(p)$ is infinite. 
\end{proof}
We leave to our gentle reader the (rather nontrivial) task
of deriving {\it{necessary}} 
conditions for $S_2(m)$ to be nonempty.
As we shall see in the next section, the behavior of $S_k(m)$ is
dramatically different when $k \ge 3$.

\section{Thue equations}

A famous theorem of the Norwegian mathematician Axel Thue \cite{Th} asserts, if
$\theta$ is an algebraic number of degree $n \geq 3$ and $\epsilon > 0$, that the inequality
$$
\left| \theta - \frac{x}{y} \right| < \frac 1{|y|^{n/2+1+\epsilon}}
$$
has at most finitely many solutions in integers $x$ and $y$ with $y
\neq 0$. Note that, if $\theta = \sqrt[k]{m}$, then we have the
algebraic identity 
$$
|x^k-my^k| = y^k \; \left| \theta - \frac{x}{y} \right| \cdot \left(
\left( \frac{x}{y} \right)^{k-1} + \theta\left( \frac{x}{y}
\right)^{k-2} + \cdots + \theta^{k-1} \right).
$$
If $m$ is not a perfect $k$-th power and $c \neq 0$ is an integer,
then it follows that the equation
\begin{equation} \label{thuebaby}
x^k-my^k = c
\end{equation}
has at most finitely many solutions in integers $x$ and $y$.
Such equations are nowadays termed
{\em Thue equations} (a nice exposition of this active area of research
may be found in the book of Fel'dman and Nesterenko \cite{FN}). 

In our situation, this immediately implies that $S_k(m)$ is finite for
each fixed integer $k \geq 3$. 
As we shall see in the next two sections, it is possible to ``effectively'' determine each such $S_k(m)$,
(and to derive an upper bound upon $k$ that depends solely upon $m$.

\section{Linear forms in logarithms}

The sets $S_k(m)$ are determined by equalities of the shape
$$
u^k - m v^k = \pm \frac{m}{d} k
$$
for $u$ and $v$ positive integers.
In all cases, we thus have
$$
\left| u^k - m v^k \right| \leq mk.
$$
It follows that
$$
\left| m^{-1} \biggl(\frac uv\biggr)^k - 1 \right| \leq \frac{k}{v^k}
$$
and hence the {\it{linear form in logarithms}}
$
\left| k \log (u/v) - \log (m) \right|
$
is ``small''. A classical result of Gel'fond \cite{Gel} (extending his work on Hilbert's
7th problem) indicates that, for any given nonzero algebraic numbers $\alpha_1, 
\alpha_2, \beta_1$ and $\beta_2$, with $\log \alpha_1$ and $\log \alpha_2$ linearly independent
over the rationals, we have
$$
\left| \beta_1 \log \alpha_1 - \beta_2 \log \alpha_2 \right| \neq 0
$$
and, moreover, provides lower bounds for such a form.
Applying state-of-the-art versions of these bounds, say those due to Laurent, Mignotte and Nesterenko \cite{LMN}, we may
conclude, as in
Theorem 3 of Mignotte \cite{Mig}, that, in our situation, $k < 600$
(if $m=2$) and, more generally, as in Theorem 2 of \cite{Mig}, that
$$
k < 10676 \log m.
$$
It follows that 
$$
\bigcup_{k=3}^{\infty} S_k (m)
$$
is finite.
In the remaining sections, we will describe a strategy for explicitly
determining this set and illustrate
it in  the cases $m=2$ and $m=3$.

\section{Solving the remaining equations}

For small values of $k$, it is possible to use standard computational techniques
based upon lower bounds for linear forms in logarithms, combined with lattice basis
reduction,
to solve the Thue equations that occur (a reasonably accessible
book which covers this field is that of Smart \cite{Sm}). If, however, $k$ is moderately large, this becomes
a computational problem due to the difficulty in finding systems of independent units in
the ring of integers of $\mathbb{Q} ( \sqrt[k]{m} )$. (If one has cycles to burn, try, for instance, to find the fundamental units
in, say, $\mathbb{Q} \left( \sqrt[41]{2} \right)$.)

In our situation, though, we are able to find local obstructions to solvability for virtually all
values of $k$ under consideration, obviating the need for extensive computations.
For the sake of simplicity, let us restrict our attention to the case $m=2$ (where, as mentioned previously, we may
assume $k < 600$). Here the equations to
be studied are of the shape
$$
x^k - 2y^k = \pm 2^{\delta} k
$$
where  $\delta \in \{ 0, 1 \}$.
Lemma \ref{lem1} also allows us to suppose that $k$ is odd.
For each such $k$, we consider primes of the form $p=2nk+1$ for $n \in \mathbb{N}$. For these $p$,
there are at most $(2n+1)^2$ values of $x^k - 2y^k$ modulo $p$. If none of these are congruent to 
$\pm k$ or $\pm 2k$ modulo $p$, we have found a local obstruction to solvability and can thus conclude that
$S_k(2)$ is empty.
If, for example, $k=13$, we consider the above equation(s) modulo $53$.
Noting that $x^{13} \equiv 0, \pm 1, \pm 23 \mod{53}$, it follows that
$$
x^{13} - 2 y^{13} \equiv 0, \pm 1, \pm 2, \pm 3, \pm 6, \pm 7, \pm 8, \pm 16, \pm 21, \pm 23, \pm 25
\mod{53}
$$
and hence may conclude
that $S_{13} (2)$ is empty.
Similar arguments suffice to eliminate all equations (for $k < 600$)
but
$$
x^3-2y^3= \pm 3, \hskip2ex x^3-2y^3= \pm 6, \hskip2ex x^5-2y^5 = \pm 10,
$$
$$
x^7-2y^7 = \pm 7, \hskip2ex x^7-2y^7 = \pm 14 \hskip2ex \mbox{ and } \hskip2ex x^{11} - 2 y^{11} = \pm 22
$$
(to verify this, the alert reader may wish to write her own piece of code;
the cognizant author used Maple \texttrademark  \cite{Maple}, being
careful to avoid our  
silico-aceric chum's rather dubious
``msolve'' routine).

Similarly, when $m=3$, we easily deal with all the equations encountered, with
the exceptions of 
$$
x^3-3y^3 = \pm 3, \hskip2ex x^5-3y^5 = \pm 15 \hskip2ex \mbox{ and } \hskip2ex x^7 - 3 y^7 = \pm 21
$$
(though we have many more $k$ to treat -- up to $k = 11728$, including even values).

\section{Endgame}

As mentioned in the last section, there are various techniques from Diophantine approximation that may be used
to solve these remaining equations. The symbolics package Kash \cite{Kash},
for example, has a built-in Thue solver
that can handle all the equations under consideration in a matter of minutes on a Sun Ultra. In any case, nowadays it is
a routine matter to verify that of the equations that have so far evaded our net,
only those with $k=3$ possess solutions, corresponding to $5^3-2 \cdot 4^3 = -3$, 
$2^3 - 2 \cdot 1^3 = 6$ and $3^3-3 \cdot 2^3 = 3$.
By way of example, to solve the 
Diophantine equations $|x^3 - 2y^3| \in \{ 3, 6 \}$, one can appeal to the inequality
$$
\left| x^3 - 2y^3 \right| \geq \sqrt{|x|},
$$
valid for all integers $x$ and $y$ (see e.g. \cite{Be}).

The equation $5^3-2 \cdot 4^3 = -3$ maps back to the solution $(x,y)$ to $x^y = y^{2x}$
given in (\ref{eq3}). Similarly, $2^3 - 2 \cdot 1^3 = 6$  yields the solution $(x,y) = (2,16)$,
also to $x^y = y^{2x}$.
The equation $3^3-3 \cdot 2^3 = 3$ which potentially yields
a solution to $x^y = y^{3x}$, 
leads to $a=27, b=24$, contradicting the coprimality of $a$ and $b$.

Similar arguments suffice to completely solve $x^y = y^{mx}$ for, with a modicum
of computation, all values of $m$ up to $40$ or so. As $m$ increases, we are faced
with the prospect of solving Thue equations of higher and higher degree, an
apparently formidable task. We note that $S_k(m)$ can be nonempty for arbitrarily
large $k$ (for example, this is always the case for $S_k (2^k-k)$ if $k$ is odd).

\section{Postscript : Mutterings on local-global principles}

A standard heuristic employed in the field of Diophantine problems is that an equation
should be solvable over 
$\mathbb{Q}$ precisely when it is solvable over $\mathbb{R}$ and over $p$-adic fields $\mathbb{Q}_p$,
for all primes $p$. Such {\it{local-global}} or {\it{Hasse}} principles are known to be true in various settings, but 
false in general.

One of the first instances where local-global principles were shown to fail
occurs in work of Skolem \cite{Sk} of 1937, where he demonstrated that the equation
$x^3-2y^3=47$ can be solved  modulo $p^{\alpha}$ for every prime power
$p^\alpha$, but has no solution over the integers.

In our analysis, we encounter numerous equations for which a like conclusion holds,
of degree up to (at least) $19$. Such is the case, for instance, for the equations
$$
x^3-5y^3=15, \hskip2ex x^5-2y^5 = 10, \hskip2ex x^5-3y^5 = 15, \hskip2ex
x^7 - 2y^7 = 7,
$$
$$
x^7-2y^7=14, \hskip2ex x^7-3y^7 = 21, \hskip2ex x^{11}-2y^{11} = 22.
$$
To actually prove that the Hasse principle fails in these instances is an interesting exercise which we
leave to our
loyal (and, at this
stage, perhaps rather fatigued)
reader!
\section{Acknowledgments}
We thank Gerry Myerson for his comments on an earlier version of this
paper. 
Some computations in this paper were performed by using Maple\texttrademark.
Maple and Maple 8 are registered trademarks of Waterloo Maple Inc.


\end{document}